\documentclass[leqno, fleqn]{amsart}%
\usepackage{amsmath}
\usepackage{amsfonts}
\usepackage{amssymb}%

\makeatletter \DeclareMathSymbol{\Gamma}{\mathalpha}{letters}{"00}
\DeclareMathSymbol{\Delta}{\mathalpha}{letters}{"01}
\DeclareMathSymbol{\Theta}{\mathalpha}{letters}{"02}
\DeclareMathSymbol{\Lambda}{\mathalpha}{letters}{"03}
\DeclareMathSymbol{\Xi}{\mathalpha}{letters}{"04}
\DeclareMathSymbol{\Pi}{\mathalpha}{letters}{"05}
\DeclareMathSymbol{\Sigma}{\mathalpha}{letters}{"06}
\DeclareMathSymbol{\Upsilon}{\mathalpha}{letters}{"07}
\DeclareMathSymbol{\Phi}{\mathalpha}{letters}{"08}
\DeclareMathSymbol{\Psi}{\mathalpha}{letters}{"09}
\DeclareMathSymbol{\Omega}{\mathalpha}{letters}{"0A}
\DeclareMathSymbol{\varGamma}{\mathalpha}{operators}{"00}
\DeclareMathSymbol{\varDelta}{\mathalpha}{operators}{"01}
\DeclareMathSymbol{\varTheta}{\mathalpha}{operators}{"02}
\DeclareMathSymbol{\varLambda}{\mathalpha}{operators}{"03}
\DeclareMathSymbol{\varXi}{\mathalpha}{operators}{"04}
\DeclareMathSymbol{\varPi}{\mathalpha}{operators}{"05}
\DeclareMathSymbol{\varSigma}{\mathalpha}{operators}{"06}
\DeclareMathSymbol{\varUpsilon}{\mathalpha}{operators}{"07}
\DeclareMathSymbol{\varPhi}{\mathalpha}{operators}{"08}
\DeclareMathSymbol{\varPsi}{\mathalpha}{operators}{"09}
\DeclareMathSymbol{\varOmega}{\mathalpha}{operators}{"0A}

\newcommand{\allmodesymb}[2]{\relax\ifmmode{\mathchoice
{\mbox{\fontsize{\tf@size}{\tf@size}#1{#2}}}
{\mbox{\fontsize{\tf@size}{\tf@size}#1{#2}}}
{\mbox{\fontsize{\sf@size}{\sf@size}#1{#2}}}
{\mbox{\fontsize{\ssf@size}{\ssf@size}#1{#2}}}} \else
\mbox{#1{#2}}\fi}

\makeatother
\makeatletter
\renewcommand*\subjclass[2][2000]{%
  \def\@subjclass{#2}%
  \@ifundefined{subjclassname@#1}{%
    \ClassWarning{\@classname}{Unknown edition (#1) of Mathematics%
      Subject Classification; using '2000'.}%
  }{%
    \@xp\let\@xp\subjclassname\csname subjclassname@#1\endcsname%
  }%
} \makeatother

\theoremstyle{plain}

\theoremstyle{remark}

\allowdisplaybreaks \numberwithin{equation}{section}
\begin{document}
\title{Lu Qi-Keng Conjectue and Hua Domain}
{\bf Dedicated to Lu Qi-Keng on the occasion of his 80th birthday}

\author{ Weiping YIN}
\address{W YIN: Dept. of Math., Capital Normal Univ., Beijing 100037, China}
\email{wyin@mail.cnu.edu.cn; wpyin@263.net} \subjclass{32H10, 32F15}
\keywords{Lu Qi-Keng conjecture, Zeroes of Bergman kernel function,
Hua domain, Bergman kernel function, K\"{a}hler-Einstein Metric,
Bergman Metric, Metric Equivalence. }
\thanks{Project supported in part by NSF of China (Grant NO.
10471097) and Specialized Research Fund for the Doctoral Program of
Higher Education of NEM of China.}

\begin{abstract}
The first part I talk about the motivation for Lu Qi-Keng conjecture
and the results about the presence or absence of zeroes of the
Bergman kernel function of a bounded domain in ${\bf{C^n}}$. The
second part I summarize the main results on Hua domains, such as the
explicit Bergman kernel function, comparison theorem for the
invariant metrics, explicit complete Einstein-K\"ahler metrics, the
equivalence between the Einstein-K\"ahler metric and the Bergman
metric etc.

\end{abstract}

\maketitle

\section*{I. Lu Qi-Keng Conjecture}
Lu Qi-Keng conjecture comes on a problem in Lu's paper "The K\"ahler
manifolds with constant curvature"[1]. That paper is written by
Chinese and published in Acta Mathematica Sinica, 1966,16:269-281.
Which was the last one before the culture revolution. Fortunately
that paper was chosen to be translated in English and published in
"Chinese Mathematics"(1966,8: 283-298). Therefore that paper was
known in the world. That problem is called Lu Qi-Keng conjecture
firstly in 1969 by M.Skwarczynski in his paper "The invariant
distance in the theory of pseudoconformal transformations and the Lu
Qi-Keng conjecture"[2]. Since then, there are many mathematicians to
research the Lu Qi-Keng conjecture, for example, Harold P.Boas, Emil
J. Straube. In 1979, R.E.Greene and H.Wu expatiated the Lu Qi-Keng
conjecture in their monograph "Function theory on manifolds which
posses a pole"(published by Springer Verlag, Lecture Notes in
Mathematics 699). 1982 S.G.Krantz expatiated the Lu Qi-Keng
conjecture in his monograph "Function Theory of Complex
Variables"(John Wiley £¦ Sons,P.53). In 1993, M.Jarnicki and P.Pflug
expatiated the Lu Qi-Keng conjecture in their monograph "Invariant
Distances and Metrics in Complex Analysis" (Walter de Gruyter,
P.184). H.P.Boas gave a lecture titled "Lu Qi-keng's Problem" in
International Conference KSCV3 in 1998, which is published in
2000[3]. In his excellent survey summarizes the various results and
methods on the Lu Qi-Keng conjecture since 1969 and six open
problems are stated. He said: "It is a difficult to determine
whether the Bergman kernel function of a specific domain has zeroes
or not. If the kernel function is presented as an infinite series,
then locating the zeroes may be of the same order of difficulty as
proving the Riemann hypothesis; and even if the series can be summed
in closed form, determining whether or not $0$ is in the range may
be hard."

Lu Qi-Keng agrees to organize an International Workshop on "Zeroes
of the Bergman kernel function" in the next year, which will be
sponsored by the Academy of Mathematics and System Sciences, Chinese
Academy of Sciences and the Capital Normal University and etc.

{\bf I.1.  Bergman kernel function}

Let $D$ be a bounded domain in ${\bf{C^n}}$, $H(D)$ be the set of
all holomorphic functions on $D$, then
$$
 L^2_a(D)=\{ f(Z)\in H(D): \int_D|f(Z)|^2dV<\infty \}
$$ is a
Hilbert space, where $dV$ is the Euclidean measure on domain $D$.
The inner product $<f,g>$ is defined by
$$
\int_D f(Z)\overline{g(Z)}dV,  f, g\in L^2_a(D).
 $$
It has countable basis. Let $\{\phi_k(Z)\}, k=1,2,\dots$ be the
complete orthonormal basis of $L^2_a(D)$, then
$$
 K_D(Z,\overline{T})=\sum_{k=1}^{\infty}\phi_k(Z)\overline{\phi_k(T)}, (Z,T)
 \in D\times D \eqno{(1.1)}
 $$
 is called the Bergman kernel function of $D$. It is unique and
 independence of the choice of complete orthonormal basis. It is
 also called reproducing kernel because of  its
 reproducing property: $\forall f\in L^2_a(D)$,
 $$
 f(Z)=\int_{D}f(T)K_D(Z,\overline{T})dV_T, Z\in D \eqno{(1.2)}
 $$
 This is a characteristic property of the Bergman kernel function
 on $D$, that is , formula (1.2) can be regarded as the definition
 of the Bergman kernel function on $D$.

 If $F: D_1 \to D_2$ is a biholomorphic mapping, and $K_1$ and $K_2$ denote the Bergman kernel functions of the Domains $D_1$
 and $D_2$ in ${\bf{C^n}}$ respectively, let $J_F(Z)$ be the
 Jacobian of $F(Z)$,  then
 $$
 K_1(Z, W)=detJ_F(Z)K_2(F(Z),F(W))det\overline{J_F(W)}.
 \eqno{(1.3)}
 $$ Where $det$ is the determinant. This relationship holds because
 if $\{\phi_j\}$ is a complete orthonormal basis of $L^2_a(D_2)$,
 then $\{(det J_F)\phi_j\circ F\}$ is a complete orthonormal basis of
 $L^2_a(D_1)$.

If $D_1=D_2=D$ and $K(Z,W)$ is the Bergman kernel function of $D$,
then the (1.3) becomes
$$
 K(Z, W)=detJ_F(Z)K(F(Z),F(W))det\overline{J_F(W)}.
 \eqno{(1.4)}
 $$
If $D$ is the homogeneous complete circular domain and $Z_0$ is any
point in $D$, then there exists a holomorphic automorphism $F$ of
$D$ such that $F(Z_0)=0.$ At this time the (1.4) becomes
$$
 K(Z_0, W)=detJ_F(Z_0)K(0,F(W))det\overline{J_F(W)}.
 $$
It is well known that if the complete orthonormal basis of Hilbert
space of square-integrable holomorphic functions on $D$ is
$\{\varphi_k(Z)\}_{k=0,1,2,\dots}$, then each $\varphi_k(Z)$ is the
homogeneous polynomial of degree $k$ in $Z=(z_1,z_2,\dots,z_n)$. And
$\varphi_0(Z)$ is the constant $a_0^1\neq 0$.

Therefore $K(0,F(W))=|a^1_0|^2=K(0,0)$. Because the $Z_0$ is the any
point in $D$, let $Z_0=Z$ in (1.4),  then one can get the formula of
Bergman kernel function for the bounded homogeneous complete
circular domain in $\bf{C^n}$:
$$
 K(Z, W)=detJ_F(Z)K(0,0)det\overline{J_F(W)}.
 \eqno{(1.5)}
 $$
Where $1/K(0,0)$ is the volume of $D$.

{\bf I.2.  Motivation for Lu Qi-Keng Conjecture}

The Riemann mapping theorem characterizes the planar domains that
are biholomorphically equivalent to the unit disk. In the higher
dimensions, there is no Riemann mapping theorem, and the following
problem arise:

Are there canonical representatives of biholomorphic equivalence
classes of domains?

In the dimension one, if $K(z,w)$ is the Bergman kernel function of
simply connected domain $D\neq {\bf{C}}$, it is well known that the
biholomorphic mapping
$$F(z)=\frac{1}{K(t,t)}\frac{\partial}{\partial \overline{w}}log\frac{K(z,w)}
{K(w,w)}|_{w=t}$$ maps the $D$ onto unit disk.

In the higher dimensions, Stefan Bergman introduced the notion of a
"representative domain" to which a given domain may be mapped by
"representative coordinates". If $D$ is a bounded domain in
${\bf{C^n}}$, $K(Z,W)$ is the Bergman kernel function of $D$, let
$$T(Z,Z)=(g_{ij})=(\frac{\partial^2logK(Z,W)}{\partial z_i\partial \overline{z_j}})$$ and its converse is
$T^{-1}(Z,W)=(g^{-1}_{ji})$. Then the local representative
coordinates based at the point $t$ is
$$f_i(Z)=\sum^n_{j=1}g^{-1}_{ji}\frac{\partial}{\partial{\overline{W_j}}}
log\frac{K(Z,W)}{K(W,W)}\mid_{W=t}, i=1,\dots,n.\eqno{(1.6)}$$ Or

$$F(Z)=(f_1,\dots,f_n)=\frac{\partial}{\partial{\overline{W}}}
log\frac{K(Z,W)}{K(W,W)}\mid_{W=t}T^{-1}(t,t),\eqno{(1.7)}$$ where
$$\frac{\partial}{\partial{\overline{W}}}=(\frac{\partial}{\partial{\overline{W_1}}},
\frac{\partial}{\partial{\overline{W_2}}}, \dots,
\frac{\partial}{\partial{\overline{W_n}}}).$$ These coordinates take
$t$ to $0$ and have complex Jacobian matrix at $t$ is equal to the
identity.

Zeroes of the Bergman kernel function $K(Z,W)$ evidently pose an
obstruction to the global definition of Bergman representative
coordinates. This observation was Lu Qi-Keng's motivation for asking
which domains have zero-free Bergman kernel functions. This problem
is called Lu Qi-Keng conjecture by M.Skwarczynski in 1969 in his
paper [2]. If the Bergman kernel function of $D$ is zero-free, that
means the Lu Qi-Keng conjecture has a positive answer, then the
domain $D$ is called the Lu Qi-Keng domain.

{\bf I.3.  Results for Lu Qi-Keng Conjecture}

The following results are well known and very easy to prove.

A. It is well known that the Bergman kernel function of a product
domain is the product of the Bergman kernel functions of the lower
dimensional domains. Therefore the Bergman kernel function of the
Cartesian product of a Lu Qi-Keng domain with a Lu Qi-Keng domain is
zero-free, and the Bergman kernel function of the Cartesian product
domain of a non Lu Qi-Keng domain with any domain does have zeroes.

B. Due to the transformation rule (1.3), the range of the
biholomorphic mapping of a Lu Qi-Keng domain is also a Lu Qi-Keng
domain.

C. If $D_j$ form an increasing sequence whose union is $D$, and
$K_m(z,w), K(z,w)$ are the Bergman kernel functions of $D_j, D$
respectively, it is well known that $lim K_j(z,w)=K(z,w)$. Due to
the Hurwitz's theorem, if the Bergman kernel function of the
limiting domain $D$ has zeroes, then so does the Bergman kernel
function of the approximating $D_j$ when $j$ is sufficiently large.

D. Due to the transformation rule (1.5), the Bergman kernel function
of any bounded homogeneous complete circular domain is zero-free.
Therefore the Bergman kernel function of the unit disk is evidently
zero-free. Consequently, the Bergman kernel function of every simply
connected planar domain$\neq \bf{C}$ is zero-free by the Riemann
mapping theorem.

{\bf I.3.1.  In the case of Bergman kernel function has zeroes}

The first counterexample for Lu Qi-Keng conjecture is given by
M.Skwarczynski in 1969,  he find  a multiply connected domain $D$ in
$\bf{C}$ such that the domain $D$ does not satisfy the Lu Qi-Keng
conjecture[2].

1969, Paul Rosenthal[12] proved that the Bergman kernel function of
annulus $D=\{0<r<|z|<1\}$ has zeroes; and for $k>2$ there exists a
bounded planar domain of connectivity $k$ which is not a Lu Qi-Keng
domain.

1976, in their paper [14], Nobuyuki Suita and Akira Yamada proved
that every finite Riemann surface which is not simply-connected is
never a Lu Qi-Keng domain.

1986, Harold P.Boas[13] proved that there exists a smooth bounded
strongly pseudoconvex complete Reihardt domain in $\bf{C^2}$, whose
Bergman kernel function has zeroes.

1998, P.Pflug and E.H.Youssfi proved [10] that the minimal
ball=$\{z\in {\bf{C^n}}: |z|^2+|zz'|<1 \}$ is not Lu Qi-Keng domain
if $n\geq 4$. Its explicit Bergman kernel function can be found in
[11]. Furthermore, when $0<a\leq 1$, and $m$ is a sufficiently large
integer, the interior approximating domain defined by
$(|z|^2+|zz'|)^m+(|z|^2-|zz'|)^m+a^m|z|^{2m}<1$ is a concrete
example of a smooth, strongly convex, algebraic domain whose Bergman
kernel function has zeroes.

1999, in their paper [16], Harold P.Boas, Siqi Fu, Emil J. Straube
proved that the Bergman kernel function of domain $\{z\in
{\bf{C^2}}: |z_1|+|z_2|^{2/p}<1\}$ has zeroes if $p>2$. In the
higher dimension they proved that the Bergman kernel function of
convex domain $\{z\in {\bf{C^n}}: |z_1|+|z_2|+\dots+|z_n|<1\}$ in
${\bf{C^n}}$ has zeros if and only if $n\geq 3$ (because if $n=2$
its Bergman kernel function is zero-free). They proved also that the
convex domain $\{z\in {\bf{C^n}}: |z_1|+|z_2|^2+|z_3|^2+
\dots+|z_n|^2<1 \}$ is not Lu Qi-Keng domain if and only if $n\geq
4$. They ask that exhibit a bounded convex domain in $\bf{C^2}$
whose Bergman kernel function has zeroes in the interior of the
domain.

2000, in his paper [5], Nguy$\hat{e}$n Vi$\hat{e}$t Anh proved that
there exists a strongly convex algebraic complete Reihardt domain
which is not Lu Qi-Keng in ${\bf{C^n}}$ for any $n \geq 3$.

2000, in his paper [9], Miroslav Engli$\check{s}$ proved that the
Hartogs domain $\{(z,t)\in \Omega \times {\bf{C^m}}: \|t\| <F(z)\}$
is not Lu Qi-Keng domain, if $m$ is a sufficiently large
integer(deponding on $F$), where $\Omega$ is a bounded, strongly
convex domain in $\bf{C^d}$ with $C^{\infty}$ smooth boundary, $-F$
is a strongly convex $C^{\infty}$ defining function for $\Omega$.
But there is no control on the size of $m$.

2000, in his survey lecture [3], Harold P.Boas point out that the
Bergman kernel function $K_{\Omega}(z,w)$ of domain
$\Omega_H=\{(z_1,z_2)\in {\bf{C^2}}: |z_2|<\frac{1}{1+|z_1|}\}$ has
a zero at origin, that is $K_{\Omega}(0,0)=0$. This fact is not true
if $\Omega_H$ is bounded.

However, due to the formula (1.5) the bounded homogeneous complete
circular domains are always the Lu Qi-Keng domain.

{\bf I.3.2. In the case of Bergman kernel function is zero-free}

1996, in his paper [15], Harold P.Boas proved that the bounded
domains of holomorphy in ${\bf{C^n}}$ whose Bergman kernel functions
are zero-free form a nowhere dense subset (with respect to a variant
of the Hausdorff distance) of all bounded domains of holomorphy.
Thus, contrary to former expectations, it is the normal situation
for the Bergman kernel function of a domain to have zeroes.

1973, in his paper [8], Shozo Matsuura proved that any bounded
complete circular domain with center at the origin is a Lu Qi-Keng
domain. 1982, in his paper [7], Tadayoshi Kanemaru proved that the
domain $D_p=\{z\in {\bf{C^n}}: |z_1|^{2/p_1}+\dots+|z_n|^{2/p_n}<1,
\}$ is a Lu Qi-Keng domain. 1979, in his paper [6], Shigeki Kakurai
use the theorem 1 in [8] proved that if $D$ is itself a
simply-connected bounded and minimal domain with center at the
origin, then $D$ is the Lu Qi-Keng domain.

In fact, above three results in [6,7,8] are not true. A domain is
called complete circular if whenever it contains a point $z$, it
also contains the one-dimensional disk $\{\lambda z: |\lambda|\leq
1\}$.

Therefore the domain $D_p=\{z\in {\bf{C^n}}:
|z_1|^{2/p_1}+\dots+|z_n|^{2/p_n}<1\}(p_j\in \mathbb{N})$ is the
bounded complete circular domain.

1999, Harold P.Boas, Siqi Fu, Emil J. Straube proved that the domain
$\{z\in {\bf{C^2}}: |z_1|+|z_2|^{2/p}<1 \}$ is not Lu Qi-Keng domain
if $p>2$, it contradicts the facts proved in [7,8]. And in 1998,
P.Pflug and E.H.Youssfi proved that the minimal ball$\{z\in
{\bf{C^n}}: |z|^2+|zz'|<1 \}$ is not Lu Qi-Keng domain if $n\geq 4$,
this negates the fact proved in [6].

It is 40 years over since the Lu Qi-Keng conjecture appeared in
1966. But there is no Chinese mathematician to research the Lu
Qi-Keng conjecture, why?

Recently we begin to study the Lu Qi-Keng conjecture, the following
are some elementary results obtained by us:

{\bf I.4.  Discuss the Lu Qi-keng problem on the Cartan-Hartogs
domain of the first type}

The Cartan-Hartogs domain of the first type is defined by
$$
Y_I(N,m,n;K)=\{W\in {\bf{C^N}},Z\in
R_I(m,n):|W|^{2K}<det(I-Z\overline{Z}^t),K>0\}.$$ Where
$R_I(m,n)=\{Z\in {\bf{C^{mn}}}: I-Z\overline{Z}^t>0 \}$ is the
Cartan domain of the first type, and $Z$ is $(m,n)$ complex matrix.
Firstly I consider the case of $N=1,m=1$. In this case, the
Cartan-Hartogs domain becomes
$$Y_I(1,1,n;K)=\{W\in {\bf{C}}, Z\in {\bf{C^n}}: |W|^{2K}+|z_1|^2+|z_2|^2+\dots+|z_n|^2<1.\}$$
Its Bergman kernel function is $$
K_I((W,Z);(\zeta,\xi))=K^{-n}\pi^{-(n+1)}F(Y)
(1-Z\overline{\xi})^{-(n+1/K)}.\\[3mm]
Y=(1-X)^{-1}, X=W\overline{\zeta}(1-Z\overline{\xi})^{-1/K}.
$$
Where $$ F(Y)=\sum_{i=0}^{n+1}b_i\Gamma(i+1)Y^{i+1}.
$$ And $b_0=0$, let
$$
P(x)= (x+1)[(x+1+Kn)(x+1+K(n-1))\dots(x+1+K)].
$$ Then the others
$b_i(i=1,2,\dots\dots,n+1)$ can be got by
$$
b_i=[P(-i-1)-\sum_{k=0}^{i-1}b_k(-1)^k\Gamma(i+1)/\Gamma(i-k+1)][(-1)^i\Gamma(i+1)]^{-1}.
\eqno{(3.1)}
$$
Recently, Liyou Zhang prove that above formula can be rewritten as
$$b_i=\sum_{j=1}^{i}\frac{(-1)^jP(-j-1)}{\Gamma(j+1)\Gamma(i-j+1)}.$$

It is well known that for the first type of Cartan-Hartogs domains
there exists the holomorphic automorphism $(W^*, Z^*)=F(W,Z)$ such
that $F(W,Z_0)=(W^*,0)$ if $(W,Z_0)\in Y_I$. Due to the
transformation rule (1.4), one has $$
K_I((W,Z);(\zeta,\xi))=(detJ_F(W,Z))|_{Z_0=Z}K_I[(W^*,0),(\zeta^*,0)](det\overline{J_F(\zeta,\xi)}).$$
Therefore the zeroes of $K_I((W,Z);(\zeta,\xi))$ are same as the
zeroes of $K_I[(W^*,0),(\zeta^*,0)]$. Let $W^*$ be the $W$, and
$\zeta^*$ be the $\zeta$, then we have
$$K_I[(W,0),(\zeta,0)]=K^{-n}\pi^{-(n+1)}F(y),
y=(1-W\overline{\zeta})^{-1}.\eqno{(3.2)}$$ Where
$$F(y)=\sum_{i=0}^{n+1}b_i\Gamma(i+1)y^{i+1}. \eqno{(3.3)}
$$ Where $b_j$ are same as that in (3.1).
For the detail see [16].

Please pay your attention to the following two facts:

1) If $(W,0)$, $(\zeta, 0)$, $(W^*,0)$ and $(\zeta^*, 0)$  belong to
$Y_I$, then their norms $ |W|,|\zeta|,|W^*|,|\zeta^*|$ are less than
1.

2) The diagonal of  Bergman kernel function is always great than
zero, that is $K_I((W,Z);(W,Z))>0.$

{\bf I.4.1. $Y_I(1,1,1;K)$ is Lu Qi-Keng domain}

At this time $Y_I(1,1,1;K)=\{W\in {\bf{C}}, Z\in {\bf{C}}:
|W|^{2K}+|Z|^2<1\}$, and the zeroes of its Bergman kernel function
$K_I[(W,Z),(\zeta,\xi)]$ are same as the zeroes of
$K_I[(W^*,0),(\zeta^*,0)]$. But
$$K_I[(W,0),(\zeta,0)]=K^{-n}\pi^{-(2)}F(y), F(y)=\sum_{i=0}^{2}b_i\Gamma(i+1)y^{i+1}$$
$$y=(1-W\overline{\zeta})^{-1}.$$
Where $b_1=K-1, b_2=1, b_0=0$, therefore
$$F(y)=(K-1)y^2+2y^3=y^3[(K-1)(1-W\overline{\zeta})+2].$$
But the zeroes of $F(y)$ are equal to
$W\overline{\zeta}=(K+1)/(K-1)$, its norm $|W\overline{\zeta}|>1$,
it is impossible.

{\bf {Therefore the Bergman kernel function of $Y_I(1,1,1;K)$ is
zero-free, that is the $Y_I(1,1,1;K)$ is Lu Qi-Keng domain.
Therefore we also prove that:

If $D\subset{\bf{C^2}}$ is a bounded pseudoconvex domain with real
analytic boundary and its holomorphic automorphism group is
noncompact, then $D$ is the Lu Qi-Keng domain, due to the E.Bedford
and S.I.Pinchuk's following theorem[61].

Theorem: }}If $D$ satisfies the above conditions, then $D$ is
biholomorphically equivalent to a domain of the form$$
E_m=\{(z_1,z_2)\in{\bf{C^2}}: |z_1|^{2m}+|z_2|^2<1\}$$ for some
positive integer $m$.

{\bf I.4.2. Discuss the domain $Y_I(1,1,2;K)$}

At this time $Y_I(1,1,2;K)=\{W\in {\bf{C}}, Z\in {\bf{C^2}}:
|W|^{2K}+|Z|^2<1\}$, and the zeroes of its Bergman kernel function
$K_I[(W,Z),(\zeta,\xi)]$ are same as the zeroes of
$K_I[(W^*,0),(\zeta^*,0)]$. But
$$K_I[(W,0),(\zeta,0)]=K^{-n}\pi^{-(2+1)}F(y), F(y)=\sum_{i=0}^{3}b_i\Gamma(i+1)y^{i+1}$$
$$y=(1-W\overline{\zeta})^{-1}.$$
Where $b_1=(K-1)(2K-1), b_2=3(K-1), b_3=1, $ therefore
$$F(y)=y^4[b_1(1-W\overline{\zeta})^2+2b_2(1-W\overline{\zeta})+6b_3].$$
Because $y^4\neq 0$, otherwise $|W\overline{\zeta}|=1$, it is
impossible. Let $$1-W\overline{\zeta}=t,
f(t)=b_1t^2+2b_2t+6b_3.\eqno{(3.4)}$$ Therefore the zeroes of the
Bergman kernel function of $Y_I(1,1,1;K)$ are equal to the zeroes of
$f(t)$. And the zeroes of $f(t)$ are
$$t=\frac{-b_2\pm \sqrt{b_2^2-6b_1b_3}}{b_1}.$$ Hence the zeroes of
$F(y)$ equal to
$$W\overline{\zeta}=\frac{b_1+b_2\pm \sqrt{b_2^2-6b_1b_3}}{b_1}.\eqno{(3.5)}$$
If the values of right side of above form are real numbers, due to
the second fact the values of the right side are in $(-1, 0)$; if
the values of right side are complex numbers, due to the first fact
the norm of the values are less than 1.

We will discuss the different cases for the $K$.

By calculations, the form (3.5) becomes
$$W\overline{\zeta}=\frac{b_1+b_2\pm \sqrt{b_2^2-6b_1b_3}}{b_1}=
\frac{2(K-1)(K+1)\pm \sqrt{-3(K-1)(K+1)}}{(K-1)(2K-1)}.
\eqno{(3.6)}$$ It must be in $(-1,0)$.

{\bf The case of $0<K<1/2$}

The value of (3.6) must be in $(-1,0)$, that is
$$-1< \frac{2(K-1)(K+1)\pm \sqrt{-3(K-1)(K+1)}}{(1-k)(1-2K)}<0.\eqno{(3.8)}$$
If $0<K<1/2$, the right side of form (3.8) always holds.

The lift side of the form (3.6) is that $$(1-K)(4K+1)<\pm
\sqrt{3(1-K^2)}.$$ It can not holds in the case of negative. If in
the case of "+", that is $$(1-K)(4K+1)<\sqrt{3(1-K^2)}.$$ That is
$$(1-K)(4K+1)^2<3(1+K).$$ That is $$16(K-1/2)^2(K+1/2)>0, $$ that is always holds.

Therefore in the case of $0<K<1/2$, the Bergman kernel function of
$Y_I(1,1,2;K)$ has zeroes, that is the $Y_I(1,1,2;K)$ is not Lu
Qi-Keng domain. The set of zeroes of $F(y)$ is
$$W\overline{\zeta}=\frac{2(K-1)(K+1)+ \sqrt{-3(K-1)(K+1)}}{(1-K)(1-2K)}.$$
Let $W\overline{\zeta}$ come back to $W^*\overline{\zeta^*}$, and
the $W^*$ and $\zeta^*$ are the image of $W$ and $\zeta$ under the
$F$ respectively. Therefore the zeroes of the Bergman kernel
function of $Y_I(1,1,2;K)$ in the case of $0<K<1/2$ are the
following
$$W\overline{\zeta}(1-Z\overline{\xi})^{-1/K}=\frac{2(K-1)(K+1)+ \sqrt{-3(K-1)(K+1)}}{(1-K)(1-2K)},$$
where $(W,Z; \zeta,\xi)\in (Y_I\times Y_I)$.

Therefore in the case of $0<K<1/2$, the domain $Y_I(1,1,2;K)$ is not
Lu Qi-Keng domain.

{\bf The case of $K=1/2$}

At that time, the $f(t)=2b_2t+6b_3=3(2-t)$, its zero is $t=2$, that
is $W\overline{\zeta}=-1$, due to the first fact, it is impossible.
Therefore in the case of $K=1/2$, the domain $Y_I(1,1,2;K)$ is the
Lu Qi-Keng domain.

{\bf The case of $1/2<K<1$}

At that time, $f(t)=b_1t^2+2b_2t+6$, because $0<|t|<2$, therefore
$$|f(t)|>6-4|b_1|-4|b_2|=6-4(1-K)(2K-1)-12(1-K)=2(2K-1)(2K+1)>0.$$
That means the $f(t)$ is zero-free. So the domain $Y_I(1,1,2;K)$ is
Lu Qi-Keng domain in the case of $1/2<K<1$. If $K=1$, then
$Y_I(1,1,2;K)$ is a ball, therefore $Y_I(1,1,2;1)$ is the Lu Qi-Keng
domain.

{\bf The case of $1<K<+\infty$}

At that time, the zeroes of $F(y)$ are
$$W\overline{\zeta}=\frac{b_1+b_2\pm \sqrt{b_2^2-6b_1b_3}}{b_1}
=\frac{2(K-1)(K+1)\pm \sqrt{-1} \sqrt{3(K-1)(K+1)}}{(K-1)(2K-1)}.$$
Its real part is
$$\frac{2K+2}{2K-1}>1.$$ It contradicts the first fact. Therefore in the case of $1<K<\infty$,
the domain $Y_I(1,1,2;K)$ is Lu Qi-Keng domain.

{\bf SUM UP:} In the case of $0<K<1/2$, $Y_I(1,1,2;K)$ is not the Lu
Qi-Keng domain; In the case of $1/2\leq K<+\infty$, $Y_I(1,1,2;K)$
is the Lu Qi-Keng domain.

{\bf I.4.3. Discuss the domain $Y_I(1,1,3;K)$}

At this time $Y_I(1,1,3;K)=\{W\in {\bf{C}}, Z\in {\bf{C^3}}:
|W|^{2K}+|Z|^2<1\}$, as the same reason of I.4.2 we only to discuss
the zeroes of $$K_I[(W,0),(\zeta,0)]=K^{-3}\pi^{-(2+1)}F(y),$$ $$
F(y)=\sum_{i=1}^{3}b_i\Gamma(i+1)y^{i+1}$$
$$y=(1-W\overline{\zeta})^{-1}.$$
Where $b_1=(K-1)(2K-1)(3K-1), b_2=(K-1)(11K-7), b_3=6(K-1), b_4=1.$
Let $t=W\overline{\zeta}$, therefore
$$F(y)=y^5[b_1(1-t)^3+2b_2(1-t)^2+6b_3(1-t)+24b_4].$$
Let $$f(t)=b_1(1-t)^3+2b_2(1-t)^2+6b_3(1-t)+24b_4,$$ By
calculations, one has
$$f(t)=(K+1)(2K+1)(3K+1)-(K-1)(K+1)(18K+11)t$$ $$+(K-1)(K+1)(18K-11)t^2-(K-1)(2K-1)(3K-1)t^3.$$
Then the zeroes of the Bergman kernel function of $Y_I(1,1,3;K)$ are
equal to the zeroes of $f(t)$.

In the case of $0<K<\frac{\sqrt{2}}{2}$:

At that time, $f(0)=(K+1)(2K+1)(3K+1)>0$, and
$f(-1)=48K(K-\frac{\sqrt{2}}{2})(K+\frac{\sqrt{2}}{2})<0$. Hence,
the Bergman kernel function of $Y_I(1,1,3;K)$ does has zeroes.

In case of $\frac{\sqrt{2}}{2}<K\leq 1$:

At that time, $|f(t)|\leq
[(K+1)(2K+1)(3K+1)-(1-K)(K+1)(18K+11)|-(1-K)(K+1)(18K+11)-(1-K)(2K-1)(3K-1)]
=2(24K^3-13K^2+12K-11)>0$. Hence, the Bergman kernel function of
$Y_I(1,1,3;K)$ is zero-free.

In the case of $K>1$:

At that time,
$f(1)=(K+1)(2K+1)(3K+1)-(K-1)(K+1)(18K+11)+(K-1)(K+1)(18K-11)-(K-1)(2K-1)(3K-1)
=24>0$, and
$f(\frac{K+1}{K-1})=\frac{2(K+1)}{(K-1)^2}[(K+1)^2(6K^2-12K+5)-(K-1)^2(6K^2+12K+5)]$.
But $[(K+1)^2(6K^2-12K+5)-(K-1)^2(6K^2+12K+5)]<0,$ therefore
$f(\frac{K+1}{K-1})<0.$ Hence $f(t)$ has a zero $t_1$ in
$(1,\frac{K+1}{K-1})$. Suppose the other two zeroes are $t_2$ and
$t_3$, because $f(t)$ is the polynomial with real coefficients of
order 3, therefore there are only the following two cases:

1) $t_2, t_3$ are conjugated;

2) $t_2, t_3$ are real numbers.

Due to relation between the roots and coefficients, one has
$t_1t_2t_3=\frac{(K+1)(2K+1)(3K+1)}{(K-1)(2K-1)(3K-1)}$, due to the
$t_1< \frac{K+1}{K-1}$, one has
$t_2t_3>\frac{(2K+1)(3K+1)}{(2K-1)(3K-1)}>1$. Therefore in the case
of 1), the norm of $t_2, t_3$ are great than 1, it is impossible. In
the case of 2), the $t_2, t_3$ can not be the negative numbers and
due to the two facts before I.4.1, the $t_2, t_3$ can not be in
[0,1). Therefore $t_2>1$ and $t_3>1$, these are impossible.
Therefore $Y_I(1,1,3;K)$ is the Lu Qi-Keng domain.

{{\bf SUM UP: In the case of $0<K<\sqrt{2}/2$, $Y_I(1,1,3;K)$ is not
the Lu Qi-Keng domain; In the case of $\sqrt{2}/2\leq K<+\infty$,
$Y_I(1,1,3;K)$ is the Lu Qi-Keng domain.

If $D\subset\bf{C^3}(\subset\bf{C^4})$ is a bounded pseudoconvex
domain with real analytic boundary, the Levi form has rank at least
$1(2)$ and its holomorphic automorphism group is noncompact, then
$D$ is the Lu Qi-Keng domain, due to the E.Bedford and S.I.Pinchuk's
following theorem[62].

Theorem:}} Let $D\subset{\bf{C^{n+1}}}$ be a bounded pseucoconvex
domain with real analytic boundary, and suppose that the Levi form
has rank at least $n-1$ at each point of the boundary. If its
holomorphic automorphism group is noncompact, then $D$ is
biholomorphically equivalent to a domain of the form$$
E_m=\{(w,z_1,z_2,\dots,z_n)\in{\bf{C^{n+1}}}:
|w|^{2m}+|z_1|^2+|z_2|^2+\dots+|z_n|^2<1\}$$ for some positive
integer $m$.

{\bf Remark A:}  Lu Qi-Keng was born on May 17th, 1927 in Guangdong
Province, CHINA. He graduated from Zhongshan University and moved to
Institute of Mathematics, Chinese Academy of Sciences in 1951
suggested by Professor Hua. 1980 he becomes the Academician of the
Chinese Academy of Sciences. He got many important results on
Several Complex Variables and Mathematical Physics:

1) Lu Qi-Keng(with L.-K.Hua)[50] established the "Theory of harmonic
functions of classical domains": Let $D=G/K$, an irreducible
classical domains. The boundary of $D$ decomposes into $r$ orbits of
$G$($r$ the rank of the symmetric space $G/K: B_1,\dots,B_r,$ with
$B_{j+1}\subset \overline{B_j}(j=1,\dots,r-1)$; each
$B_j(j=1,\dots,r-1)$ is a direct product space $D_j\times L_j$,
where $D_j=G_j/K_j$ is a classical symmetric space of $K$). $B_r$ is
a homogeneous space of $K$, and is the $\breve{S}$ilov boundary of
$D$. Now let there be given a real-valued continuous function $\phi$
on the $\breve{S}$ilov boundary $B_r$ of $D$, then the Dirichlet
problem for $B_r$ can be solved in the sense that there exists an
unique function $f$ in $D$ satisfying the Laplace-Beltrami equation
for the Bergman metric of $D$ and having the continuous boundary
values $\phi$ on $B_r$. It turns out that it always  has continuous
boundary values everywhere, and on each $B_j=D_j\times
L_j(j=1,\dots,r-1)$ its restriction to any cross-section above $D_j$
projects into a harmonic function on $D_j$. The proof is based on a
study of the Poisson formulas.

2) Lu Qi-Keng-Schwarz lemma on the bounded domain $D$ in
${\bf{C^n}}$[51,52]: Let $w=f(z)$ be an analytic mapping of $D$ into
${\bf{C^n}}$ such that $\sum |f_j(z)|^2\leq M^2$, then $(\partial
f/\partial z)\overline{(\partial f/\partial z)}'\leq M^2T(z,z).$
Where $M$ is a constant, and $T(z,z)$ is the Bergman metric matrix.
If $D$ is homogeneous, $F(z)$ the holomorphic mapping of $D$ into
$D$, $ds$ is the Bergman metric of $D$, then there is a constant
$k(D)$ such that $F^*ds^2\leq k^2(D)ds^2$; and $k^2(D)$ is an
analytic invariant. Let $k_1(D), k_2(D)$ be the infimum and supremum
of the holomorphic sectional curvature of $D$ under the Bergman
metric respectively, then $k^2(D)=\frac{k_1(D)}{ k_2(D)}$; therefore
one has $F^*ds^2\leq \frac{k_1(D)}{ k_2(D)}ds^2$. To the best of my
knowledge it is the first time that the Schwarz lemma contacts with
the curvatures in the higher dimension. And Lu proved that if $D$
and $G$ are the irreducible classical domains, then $D$
biholomorphic equivalent to $G$ if and only if $k(D)=k(G),
k_1(D)=k_1(G)$, or $k_1(D)=k_1(G), k_2(D)=k_2(G)$.

3) Lu Qi-Keng proved that if $D$ is the bounded domain in
${\bf{C^n}}$, then its Bergman metric always not less than its
Carath$\acute{e}$ory metric, which is imply in his paper[52].

4) In 1957 Lu Qi-Keng(with Tong-De Zhong) studied the boundary value
of the Bochner-Martinelli integral and obtained a Plemelj
formula[53].

5) Lu Qi-Keng shows that if $D$ is a bounded domain in ${\bf{C^n}}$
with a complete Bergman metric and constant holomorphic sectional
curvature, then $D$ is holomorphically isomorphic to the unit Ball
in ${\bf{C^n}}$ and the constant is equal to $-\frac{2}{n+1}$. Where
the $D$ is not need the simply connected domain. In the same time he
ask a question, which is called the Lu Qi-Keng conjecture[1].

6) Lu Qi-Keng constructs a so-called horo-hypercircle coordinate on
symmetric spaces. By using these coordinates, he is able to
establish an explicit formula for the heat kernels on the symmetric
spaces[55].

7) Lu Qi-Keng constructs an explicit global solution of the
Einstein-Yang-Mills equation on the Dirac conformal space[56] and
related topics[60].

8) Lu Qi-Keng[57] constructs a Leray map for a piecewise smooth
compact component $D$ of a domain defined by $U=\{z\in {\bf{C^N}}:
q(z,\overline{z})=0\}$ where $q(w,\overline{z})$ is a real
polynomial in $w=(w_1,\dots,w_N)$ and
$\overline{z}=(\overline{z_1},\dots,\overline{z_N})$, which does not
vanish for $(w,z)\in (D\times \overline{D})$. And he gives the
polynomials $q(w,\overline{z})$ for the classical domains and
consequently obtains an integral representation of solutions of the
$\overline{\partial}$-equation on the classical domains.

{\bf Remark B:} Besides the Lu Qi-Keng's many results above, the
foremost one is that Lu Qi-Keng[54] establishes the relation between
the theory of Yang-Mills fields and that of connections of principal
bundles. He proved that the new definition of the Yang-Mills field
suggested by Yang is equivalent to the parallelism along a curve.
Which is published in the "Acta Phys.Sinica" 1974,23(4): 249-263.
One year later C.N.Yang [59] discover the same thing and published
in the Phys.Rev.D(3), 1975,12(12):3845-3857. And in 1972, Lu Qi-Keng
wrote a Lecture Notes about this subject and gave a series lectures,
that Lecture Notes is kept by many Chinese mathematicians up to now.

\section*{II. Hua Domain}

{\bf II.1. Cartan domains.}

The Hua domains are introduced by Weiping YIN since 1998, which are
built on the Cartan domains(classical domains). Due to E.Cartan's
results,
 the irreducible bounded symmetric domains in $\bf{C}^n$
 can be categorized into four big types
 and two exceptional types.
 The four big types can be written as:$$
\begin{array}{ll}
R_I(m,n)&:=\{Z\in {\bf{C}^{mn}}: I-Z\overline{Z}^t>0, \},\\[3mm] R_{II}(p)&:=\{Z\in
{\bf{C}^{p(p+1)/2}}: I-Z\overline{Z}^t>0, \},\\[3mm] R_{III}(q)&:=\{Z\in
{\bf{C}^{q(q-1)/2}}: I-Z\overline{Z}^t>0,\},\\[3mm] R_{IV}(n)&:=\{Z\in {\bf{C}^n}:
1+|ZZ^t|^2-2Z\overline{Z}^t>0,\\& 1-|ZZ^t|^2>0\}.
\end{array}
$$
where $Z$ is $m\times n$ matrix, $p$ degree symmetric matrix, $q$
degree skew symmetric matrix and $n$ dimensional complex vector
respectively. And the two exceptional types are of dimensions 16 and
27 respectively, their matrix representations can be found in [Yin
Weiping. Two problems on Cartan domains. J of China Univ of Sci and
Tech. 1986, 16(2): 130-146.]. If $m=1$, $R_I(1,n)$ is the unit ball
in $\mathbb{C}^n$, i.e. $R_I(1,n)=B_n(0,1)=B_n$. The set of Cartan
domains is denoted by $\bf{R}_{H}$, that is
$${\bf{R}_{H}}=\{R_I(m,n), R_{II}(p), R_{III}(q), R_{IV}(n),
$$ $$R_V(16), R_{VI}(27)\}.$$

{\bf II.2. Story of Hua domains.}

While YIN were visiting the Institut Des Hautes Etudes
Scientifiques(IHES) in February 1998, Professor Guy ROOS invited YIN
to visit the University of Poitiers and gave a talk there. The title
of his talk was "Some Results on Egg Domain", the Egg domain is the
following form:
$$E(K)=\{w\in {\bf{C}}, z\in {\bf{C}^n}: |w|^{2K}+|z|^2<1\}.$$ Roos asked a
question: the Egg domain is equivalent to
$$\{w\in {\bf{C}}, z\in {\bf{C}^n}: |w|^{2K}<1-|z|^2\}.$$
The right side is $1-|z|^2>0$, that is the unit ball $B_n$,
therefore the Egg domain is equivalent to $$\{w\in \mathbb{C}, z\in
B_n: |w|^{2K}<1-|z|^2\}.$$ That is to say the Egg domain is built on
$B_n$, and $B_n$ is the special case of Cartan domains. Then Roos'
question is: can we construct a domain built on the Cartan domain?
YIN gave a positive answer as follows:
$$
Y_I(N,m,n;K)=\{W\in {\bf{C}^N},Z\in R_I(m,n):$$
$$|W|^{2K}<\det(I-Z\overline{Z}^t),K>0\},
$$ $$ Y_{II}(N,p;K)=\{W\in {\bf{C}^N},Z\in
R_{II}(p):$$ $$|W|^{2K}<\det(I-Z\overline{Z}^t),K>0\}, $$ $$
Y_{III}(N,q;K)=\{W\in {\bf{C}^N},Z\in R_{III}(q):$$
$$|W|^{2K}<\det(I-Z\overline{Z}^t),K>0\}, $$
$$Y_{IV}(N,n;K)=\{W\in {\bf{C}^N},Z\in
R_{IV}(n):$$ $$|W|^{2K}<1-2Z\overline{Z}^t+|ZZ^t|^2, K>0\},
$$
where  $\det$ indicates the determinant, $N, m, n, p, q$ are natural
numbers. These domains are called Cartan-Hartogs domains or
Super-Cartan domains. When YIN came back to Beijing, China, he
computed the Bergman kernel function in explicit formulas for all of
the Cartan-Hartogs domains, for example, the Bergman kernel function
of $Y_I(1,m,n;K)$ is
$$
K^{-mn}\pi^{-(mn+1)}F(Y) \det(I-Z\overline{Z}^t)^{-(m+n+1/K)},
$$where
$$
F(Y)=\sum_{i=0}^{mn+1}a_i\Gamma(i+1)Y^{i+1}. $$ $$Y=(1-X)^{-1},
X=|W|^2[\det(I-Z\overline{Z}^t)]^{-1/K}.
$$and $a_i$ are constants.

The right hand of above inequalities are called generic norm of
corresponding in Jordan Triple System. Let the generic norm of $D$
is denoted by $N_D(Z,\overline Z)$. Then the definition of
Cartan-Hartogs domain can be written as follows:
$$
Y(N,D;K)=\{W\in {\bf{C}^N}, Z\in D: |W|^{2K}< N_D(Z,\overline Z),
D\in \bf{R}_{H}\}.$$ After that, YIN generalizes the Cartan-Hartogs
domains to the following forms, which are called Cartan-Egg domains:
$$
\{W_1\in {\bf{C}^{N_1}}, W_2\in {\bf{C}^{N_2}}, Z\in D:
|W_1|^{2K}+|W_2|^2< N_D(Z,\overline Z), D\in \bf{R}_{H}\}.$$ YIN
generalizes the above domains to the following forms, which are
called Hua domains:
$$
\begin{array}{lll}
\{ W_j\in {\bf{C}^{N_j}},Z\in
D:\displaystyle\sum_{j=1}^r||W_j||^{2p_j}\\
< N_D(Z,\overline Z), p_j>0, j=1,\dots,r, D\in \bf{R}_{H}
\},\end{array}$$ where
$||W_j||^2=\displaystyle\sum_{k=1}^{N_j}|w_{jk}|^2$.
 The Hua domain is equivalent to the following domain:
$$
\begin{array}{lll}
\{ W_j\in {\bf{C}^{N_j}},Z\in
D:\displaystyle\sum_{j=1}^r\frac{||W_j||^{2p_j}}{N_D(Z,\overline
Z)}<1, \\p_j>0, j=1,\dots,r, D\in \bf{R}_{H}  \}.\end{array}$$ Here
all of the denominators are the same, if each denominator is
different, then one can get the following domain:
$$
\begin{array}{lll}
\{ W_j\in {\bf{C}^{N_j}},Z\in
D:\displaystyle\sum_{j=1}^r\frac{||W_j||^{2p_j}}{[N_D(Z,\overline
Z)]^{K_j}}<1, \\p_j>0, K_j>0, j=1,\dots,r, D\in \bf{R}_{H}
\}.\end{array}$$ That domain is called Hua-construction.

The Bergman kernel functions of all of Hua-construction can be
computed in explicit formula. In the following 4 cases, the Bergman
kernel functions can be expressed by the elementary functions:

(I) $r=1, p_r=p_1=K;$

(II) $p_1=p_2=\dots=p_{r-1}=1, p_r=K>0, $;

(III) $1/p_1, \dots, 1/p_r$ are positive integers;

(IV) $1/p_1, \dots, 1/p_{r-1}$ are positive integers,  $p_r>0$.

{\bf II.3. Main results on Hua domains.}

Up to now, we get some main results on Hua domains or
Hua-construction:

(I) We get the Bergman kernel functions in explicit formulas on
Cartan-Hartogs domains, Cartan-Egg domains, Hua domains and
Hua-constructions[16-37]. Our method is quite different with the
well known methods. We need the holomorphic automorphism groups of
these domains, the complete orthonormal system of Semi-Reinhardt
domain and to compute the integral $$
\int_{D}det(I-Z\overline{Z}^t)^{\lambda}dZ,$$ where $D$ is the
Cartan domain of the first type or second type or the third type. It
is well known that there are few types of domain for which the
Bergman kernel function is explicitly computable. Any domain for
which the Bergman kernel can be computed in explicit formula is a
good domain of worth researching.

(II) We get the complete Einstein-K\"{a}hler metric(or
K\"{a}hler-Einstein metric) with explicit formulas for some
Cartan-Hartogs domains(special case of Hua domains), which are
non-homogeneous and the sharp estimate of the holomorphic sectional
curvature under this metric is also given[38-43]. It was proven by
Cheng-Yau and Mok-Yau that there is an unique complete
Einstein-K\"ahler metric with Ricci curvature $-1$ on any bounded
domain of holomorphy. However, the existence of such a metric is
highly nonconstructive, and it is very difficult to write down the
metric except for homogeneous domains. As far as we are aware, the
explicit formula for Cartan-Hartogs domains are the first such
non-homogeneous examples.

(III) We proved the comparison theorem for Bergman metric and
Kobayashi metric on Cartan-Hartogs domain[44-47]. Let $D$ be the
bounded domain in ${\bf{C^n}}$, and $B_D, C_D, E_D, K_D$ be the
Bergman, Carath$\acute{e}$ory, Einstein-K\"ahler, Kobayashi metrics
respectively. It is well known that $C_D\leq K_D, C_D\leq 2B_D.$
When $D$ is any of the Cartan-Hartogs domains, we proved that the
inequalities $B_D\leq cK_D$ and $E_D\leq cK_D$ hold, where $c$ is a
constant related to $D$. The proofs of the two inequalities are
based on the following theorem [58].

{\bf Theorem.} Let $\beta$ denote an infinitesimal Finsler metric on
a $\bf{K}$-hyperbolic Banach manifold $\bf{D}$. If the holomorphic
sectional curvature of $\beta$ is bounded above by $-c^2, c\in R^+$,
then $\beta(p,v)\leq \frac{2}{c}\bf{K}_{\bf{D}}$ $(p,v)$ for all
$(p,v)\in T(\bf{D})$.

Due to the theorem, the key step is to prove that the holomorphic
sectional curvature of $B_D$ and $E_D$ are bounded above by the
negative constant($D$ is any of the Cartan-Hartogs domains). In
general it is very complicated.

(IV) We proved that the Bergman metric is equivalent to the
Einstein-K\"ahler metric on Cartan-Hartogs domains, which will be
published somewhere. Firstly, we introduce a class of new invariant
metrics on these domains, and proved that new these metrics are
equivalent to the Bergman metric, hence the new metrics are
complete. Secondly, the Ricci curvatures under these new metrics are
bounded from above and below by the negative constants. Thirdly, we
estimate the holomorphic sectional curvatures of the new metrics,
and show that they are bounded from above and below by the negative
constants. Finally, by using these new metrics as a bridge and Yau's
Schwarz lemma[48] we prove that the Bergman metric is equivalent to
the Einstein-K\"ahler metric. That means the Yau's conjecture[49] is
also true on Cartan-Hartogs domains.

{\bf Remark C:} From above we know that the Bergman kernel function
of Hua domains are computed in explicit formulas, based on these one
can research the  Lu Qi-Keng conjecture on them, that is to discuss
the the presence or absence of zeroes of the Bergman kernel function
of the Hua domains. It can be reduced to discuss the zeroes of
polynomial  on the unit disk in $\bf{C}$, or discuss the zeroes of
the convergence series $\sum_{j=1}^{\infty} b_j z^j$ on the unit
disk in $\bf{C}$, where $b_j(j=1,2, 3,\dots)$ are real numbers.

\end{document}